\documentclass[reqno,12pt]{amsart}
\usepackage{eurosym}
\usepackage{amssymb}
\usepackage{amsfonts}
\usepackage{amsmath}
\usepackage{color}
\usepackage{amsaddr}
\usepackage{graphicx}

\usepackage{caption,subcaption}

\theoremstyle{remark}

\numberwithin{equation}{section}

\begin{document}
\title[Modeling public sentiment with nonlocal diffusion]{Modeling change in public sentiment with nonlocal reaction-diffusion equations}
\author{Joseph L. Shomberg$^{*}$}
\address{$^{*}$Department of Mathematics and Computer Science, Providence
College, Providence, Rhode Island 02918, USA, \\
\texttt{jshomber@providence.edu}}
\date{\today}

\begin{abstract}
This is a brief ``proof of concept'' article that shows a nonlocal variant of a reaction-diffusion equation, which is already well suited to the study of pattern formation, is a plausible tool for the particular application of modeling change in public sentiment. 
Public sentiment is ubiquitous in modern society; from the marketing of consumer goods to the choice of political rhetoric, we are all affected by it.
Modeling change in public sentiment has an important role when analyzing and predicting the course society is on.
Of course, change is permanent. 
Perhaps the best-known feature in public sentiment is polarization.
By this, we mean the development or establishment of a certain subset of a population whose views on a given subject are fixed to one end of the spectrum.
It is when public sentiment undergoes a change do we witness the emergence of polarization.
In large part we see polarization arise through the use of persuasion and rhetoric.
Moreover, some marketing or arguments are more effective on certain individuals, thus leading these people into a polarized regime. 
Our model captures this emergence of the polarized regime. 
In addition, in some cases, we see the development of true polar opposites (which we call mixed polarity).
Our method features:
\begin{itemize}
\item We use a nonlocal Chafee--Infante reaction-diffusion equation to model the evolution of public sentiment in a population that interacts with other individuals.
\item We employ a pseudo-random convolution kernel as a symmetric matrix of lognormally distributed values. This kernel models the influence of individuals when interacting with others.
\item Change in sentiment emerges and may converge to a polarized state expressed by a double-well potential. Other more complicated states occur whereby a mixed polarization emerges. 
\end{itemize}
\end{abstract}

\maketitle

\section{Introduction and motivation}

In this article, we introduce the notion of nonlocal diffusion as a means to describe the change, and possible polarization, that occurs in some public sentiment. 
Public sentiment can be collected by a variety of sources ranging from websites that perform online surveys to Bayesian classifiers digesting our social media behavior.
Undoubtably, such data is important to advertisers, to name just one group, and the amount of data that is accessible today is increasing at an ever faster rate (cf. e.g. \cite{BPM,CFCH,CHEN2019402}).
According to a recent story at {\tt{blacklinko.com}} \cite{blacklinko-1}, ``The social media growth rate since 2015 is an average of 12.5\% year-over-year.'' 
A quick Google search will bring up thousands of pages pertaining to `modeling public sentiment', but here we are interested in its evolution brought forward through nonlocal diffusion.

Differential equations have long been used to predict the evolution of an initial state.
The size of a population with a limited amount of resources can be modeled with a logistic equation, for example.
It is sometimes important to incorporate mixing, or diffusion, into a population.
This is evident whenever a substance of high concentration prefers a region of lower concentration.
We find \cite{Volpert} analyzes a problem involving competing species where each ``consumes resources in some area around their average position''.
This later class of diffusion mechanism is now referred to as {\em nonlocal} diffusion.
One important feature in nonlocal diffusion is the interaction {\em kernel}.
Such kernels are functions where one can assign relevant interaction values, according to a certain probability distribution perhaps. 
Besides population dynamics, nonlocal diffusion has already earned its place in modeling several phenomena important to the social sciences, especially economics.
Notably, recent applications to wealth distribution appear in \cite{Anita} and \cite{Volpert2}.

Whereas (local) diffusion is championed for its success in modeling classical phenomena, such as heat flow (see J. Fourier's ``The Analytical Theory of Heat''), nonlocal diffusion is far-reaching. 
Rather than one particle interacting with only strictly adjacent particles, nonlocal diffusion provides a mechanism in which one particle can interact with all others.
The underlying rule is that one particle is allowed to interact with all others according to the probability/interaction kernel.
We will use (positive) numbers to represent the number of interactions between two individuals.
We can assume such an interaction also represents an exchange of information.
Of course, the form of that exchange could be persuasive. 
Thus, nonlocal diffusion is an ideal mathematical tool for modeling the type of interaction that often occurs when studying sentiment. 

\section{The model problem}

\begin{figure}[htpb]
\includegraphics[scale=0.76]{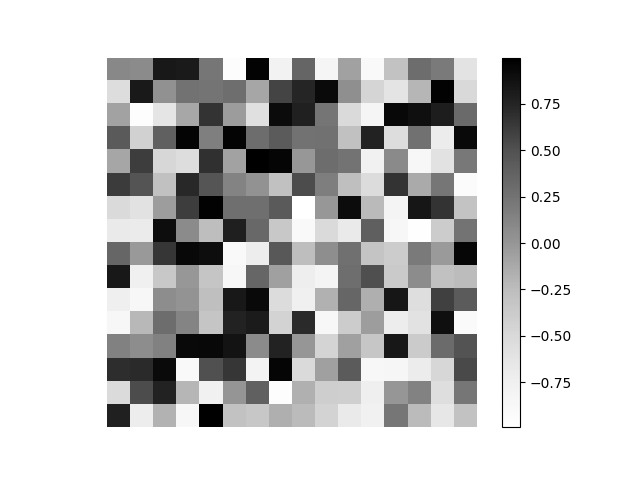}
\caption{A sample of sixteen individuals, represented as rows, and their `sentiment' to sixteen topics, in columns.}
\label{fig-1}
\end{figure}

We will suppose that sixteen people are given a questionnaire containing sixteen questions asking them for their opinion on different topics (the number of people does not have to be the same as the number of questions, but it makes for a nicer presentation when they are).
Each question is indicated by one of the sixteen (vertical) columns of Figure \ref{fig-1}.
The sixteen different people are represented as one of the sixteen (horizontal) rows.
So, the fifth row indicates the answer (or `sentiment') to each of the questions corresponding to the fifth person.
Responses can range in the interval of numbers $[-1,1]$, where the value $-1$ corresponds to ``strongly disagree'' (in this case the color white) to the value $+1$ which corresponds to ``strongly agree'' (color black).
Figure \ref{fig-1} contains sample responses.
Here, the responses are uniformly randomly chosen. 

People are allowed to discuss the topics with others who participated in the questionnaire, as well as with people outside the group of sixteen.
In this experiment, we do limit the total number of people to 31 (the sixteen who answered the same questions and fifteen who did not--the number 15 could be increased or decreased).
We assume that each individual may interact with each of the others on a certain fixed timescale (e.g. daily) and that any interaction produces a degree certain influence.
We assume {\em interaction} occurs with an exchange of communication during physical or virtual meetings.
All interactions are possible and given a nonnegative value representing the (average) number of interactions over a specified time period.
Here, {\em influence} is a certain persuasive intent.
When it comes to interacting with different people, we may barely see someone in a month, or that person could make poor arguments about certain topics and is unable to sway many others to their point of view.
On the other hand, some are easily persuaded and some seek out others in an attempt to change someone's point of view.
For the purposes of modeling change in sentiment, we do not assume one interacts with oneself.
So here, one does not change their mind on their own (however, one could easily assign a nonzero probability for this occurrence to any of the individuals).
Finally, interactions are assumed to be {\em symmetric}, which means the influence of Person A on Person B is not different from the influence of Person B on Person A.
(This assumption could also be dropped.)

\begin{figure}[htpb]
\includegraphics[scale=0.76]{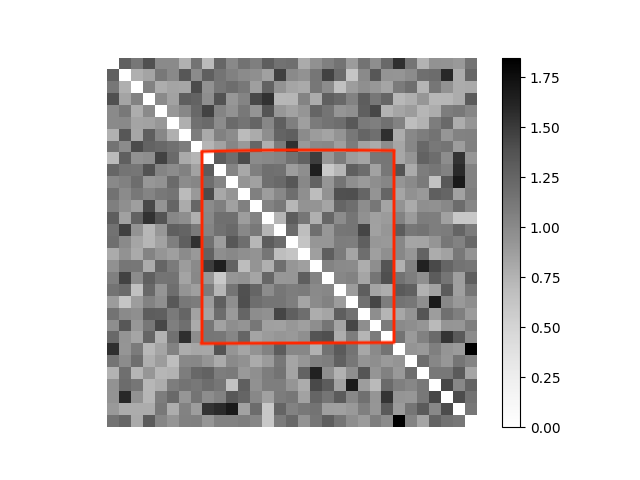}
\caption{Thirty-one individuals in rows and their interaction values with the sixteen individuals participating in the survey  (contained in the red square). The white diagonal is the zeros for self interactions.}
\label{fig-2}
\end{figure}

We need to record the influence of all the possible interactions.
For this we will use a symmetric matrix, or when displayed graphically, a $31\times 31$ array where each entry contains a number represented by some small grayscale square.
Each interaction score is a random number determined by a lognormal distribution with mean $\mu=1.0$ and standard variation $\sigma=1.7$.
This means all interactions have a nonnegative value, but increasing interaction scores occur with less frequency.
We take these interaction values to be lognormally distributed.
(It is interesting to note that the lognormal distribution is also used in pricing stock options.)
Anyway, the group of sixteen who answered the survey make up a $16\times16$ square inside the $31\times31$ interaction array shown in Figure \ref{fig-2}.

To capture polarization in public opinion, we will rely on a function that has two attractive fixed states.
The function is called a {\em double-well potential} and is given by $F(x)=\frac{1}{4}(x^2-1)^2$.

\begin{figure}[htpb]
\includegraphics[scale=0.5]{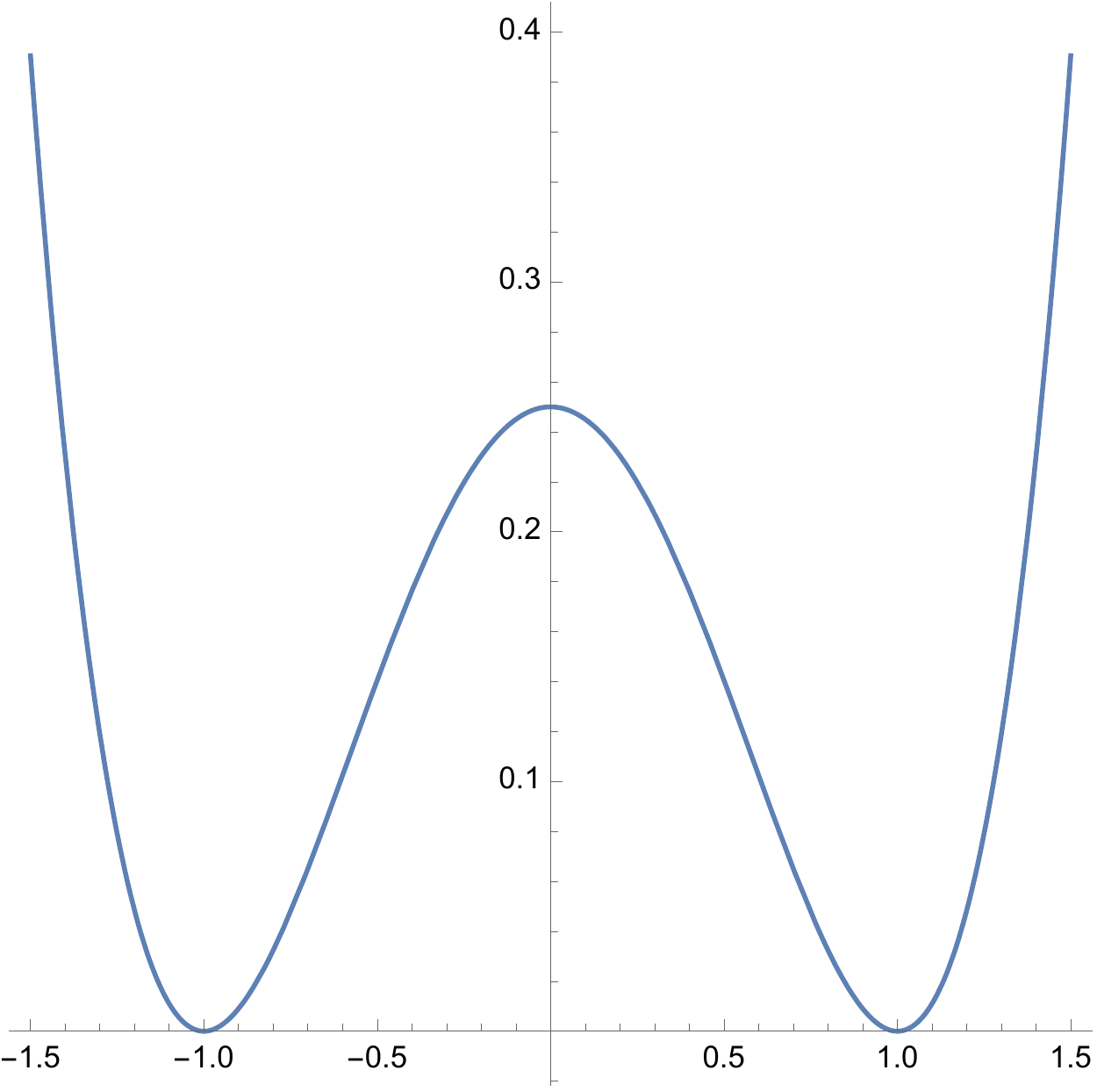}
\caption{A graph of the double-well potential $F(x)=\frac{1}{4}(x^2-1)^2$.}
\label{fig-c}
\end{figure}

Observe, $F(0)=0$ and $F(\pm1)=0.$
The numbers $x=\pm1$ refer to the {\em stable steady states} whereas $x=0$ is said to be {\em unstable}.
The graph of $F(x)$ is like a `W' and one can imagine the middle hump, $x=0$, separates sentiment into one of the two lower `wells', $x=-1$ or $x=1$.
Hence, this function works to promote positive and negative polarity.
The derivation used below actually requires the derivative of the double-well potential, $F'$ which we will denote by $f$.

In order to communicate the evolution, or change, of any given initial sentiment, such as that given in Figure \ref{fig-1} for example, we rely on a differential equation made from balancing the rate of change in time to the interaction of mixing influence in the presence of two stable steady states.
We put these concepts together by borrowing a model known to applied mathematicians as a nonlocal Chafee--Infante reaction-diffusion equation.
This equation has a celebrated history because of its archetypical impact on the development of dynamical systems (indeed, see \cite{CI}).
That is,
\begin{equation*}
\partial_tp_i + \sum_{m=1}^{16} K_m*p_i + f(p_i) = 0, \quad \text{for}\ i=1,2,\dots,16,
\end{equation*}
where $p_i=p_i(t,x)$ is the $i$-th individual's unknown sentiment to question/topic $x$ at time $t$.
The term $K_m$ represents the $m$-th column of the extended interaction kernel $K$.
The chosen convolution expresses a Neumann-type no flux boundary condition, 
\begin{equation*}
(K_m*p_i)(t,x)=\int_\Omega K_m(x-y)\left( p_i(t,y)-p_i(t,x) \right) dy.
\end{equation*}
OK, now we answer the question, ``{\em why the extra 15 people and why use the extended interaction kernel?}''
Inside the kernels $K_m$, we see the argument contains the difference $x-y$.
The way a convolution works is like one function, say the kernel, slides from left to right over another function, in this case let us simply say $p_i(t,y)-p_i(t,x)$, and we sum the area of the overlap.
So to compute the integral at $x=1$, we need to account for all $y$, i.e., $y=1,2,\dots,16.$
Hence, we need to know the value of $K_m(1-16)=K_m(-15)$ and other terms like this.
When we are at a point $x$ inside the red box in Figure \ref{fig-2}, we need to look up 15 points to find $K_m(-15)$.
Note that the boundary is ``wrapped'' in the sense that when we go out the top, we reappear at the bottom.
Moving on, in this setting the derivative of the double-well potential is
\begin{equation*}
f(p_i)=p_i^3-p_i
\end{equation*}
and the stationary polar sentiments $p_i=\pm1$ are apparent.

\begin{figure}[htpb]
    \centering 
\begin{subfigure}{0.33\textwidth}
  \includegraphics[width=\linewidth]{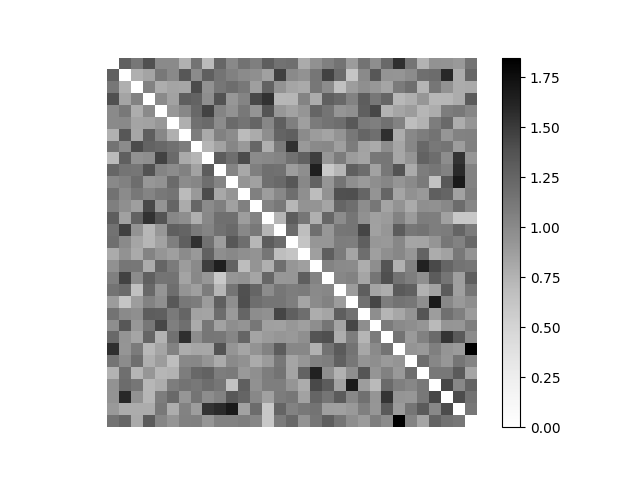}
  \caption{Extended interaction kernel}
\end{subfigure}\hfil 
\begin{subfigure}{0.33\textwidth}
  \includegraphics[width=\linewidth]{phi00000-20}
  \caption{Initial sentiment}
\end{subfigure}\hfil 
\begin{subfigure}{0.33\textwidth}
  \includegraphics[width=\linewidth]{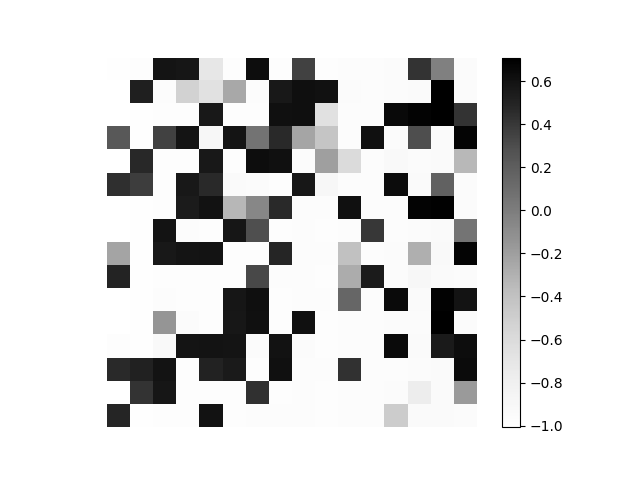}
  \caption{Sentiment at 475 iterations}
\end{subfigure}
\medskip
\begin{subfigure}{0.33\textwidth}
  \includegraphics[width=\linewidth]{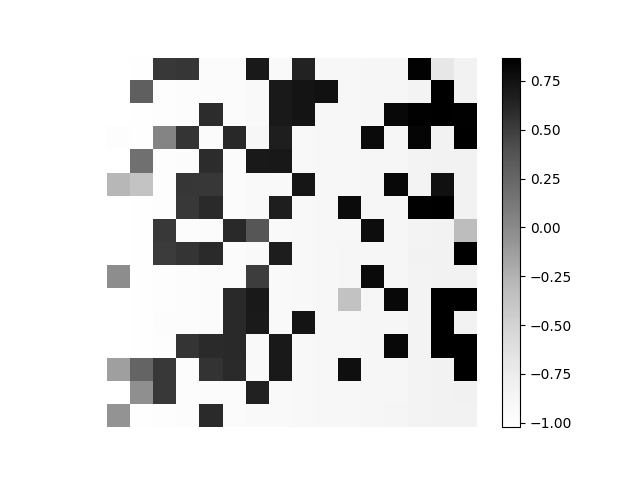}
  \caption{Sentiment at 950 iterations}
\end{subfigure}\hfil 
\begin{subfigure}{0.33\textwidth}
  \includegraphics[width=\linewidth]{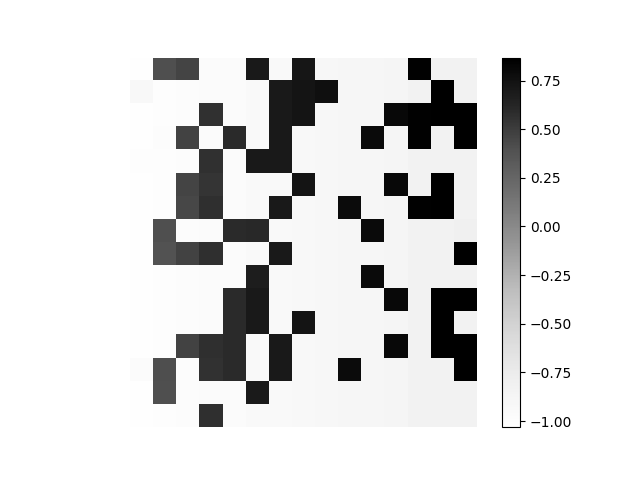}
  \caption{Final sentiment}
\end{subfigure}\hfil 
\begin{subfigure}{0.33\textwidth}
  \includegraphics[width=\linewidth]{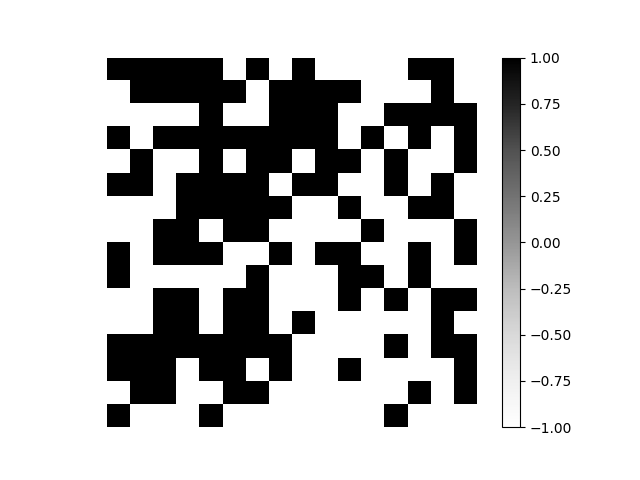}
  \caption{End state with no interaction}
\end{subfigure}
\caption{Sum of initial sentiment is $4.017$ vs. sum of final sentiment after 1425 iterations is $-121.965$. This indicates a dominantly negative (white) change in public sentiment.}
\label{fig-3}
\end{figure}

\begin{figure}[htpb]
    \centering 
\begin{subfigure}{0.33\textwidth}
  \includegraphics[width=\linewidth]{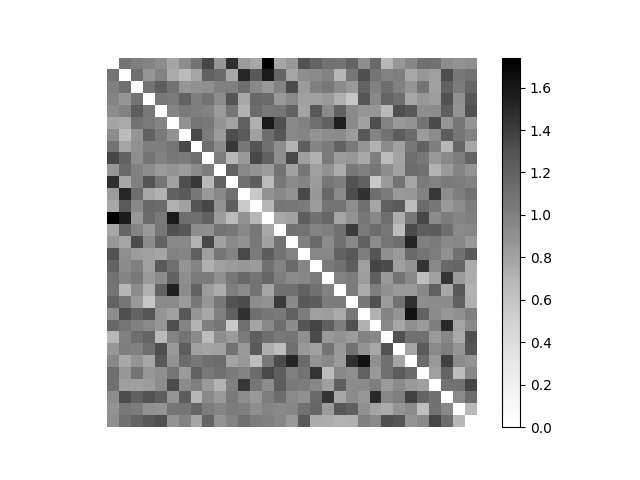}
  \caption{Extended interaction kernel}
\end{subfigure}\hfil 
\begin{subfigure}{0.33\textwidth}
  \includegraphics[width=\linewidth]{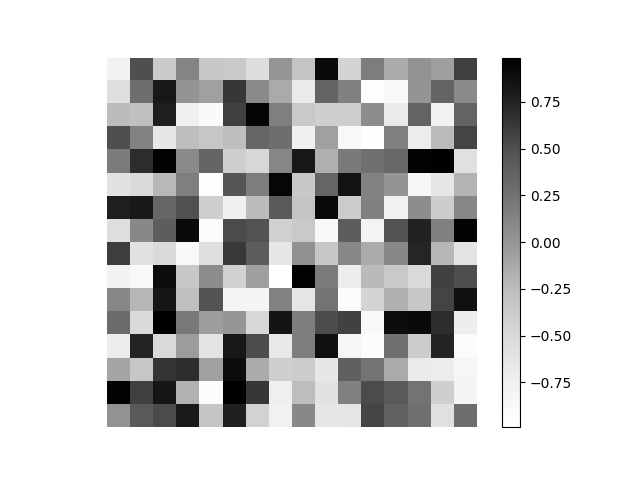}
  \caption{Initial sentiment}
\end{subfigure}\hfil 
\begin{subfigure}{0.33\textwidth}
  \includegraphics[width=\linewidth]{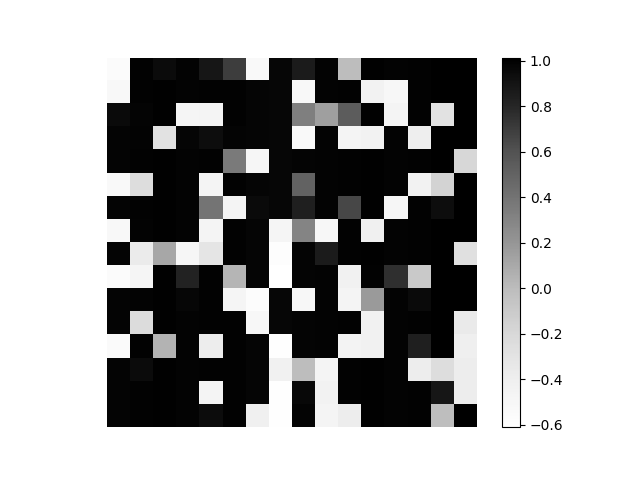}
  \caption{Sentiment at 475 iterations}
\end{subfigure}
\medskip
\begin{subfigure}{0.33\textwidth}
  \includegraphics[width=\linewidth]{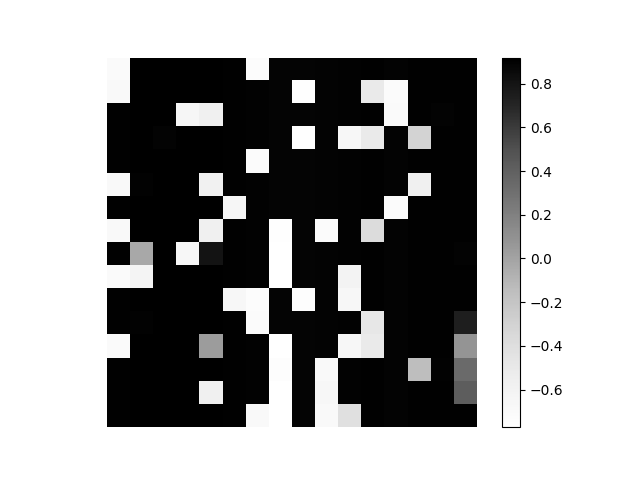}
  \caption{Sentiment at 950 iterations}
\end{subfigure}\hfil 
\begin{subfigure}{0.33\textwidth}
  \includegraphics[width=\linewidth]{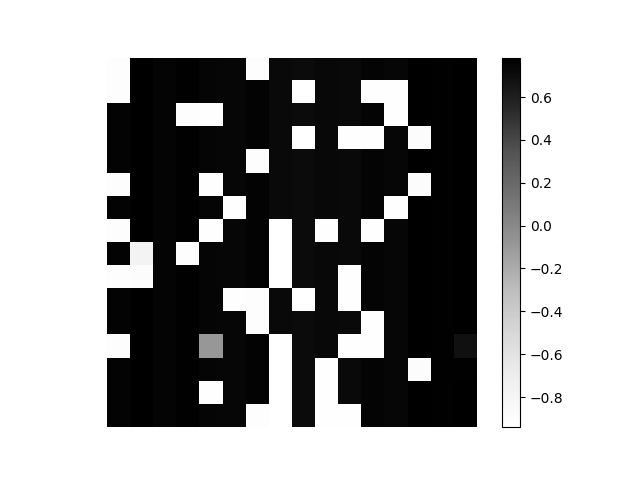}
  \caption{Final sentiment}
\end{subfigure}\hfil 
\begin{subfigure}{0.33\textwidth}
  \includegraphics[width=\linewidth]{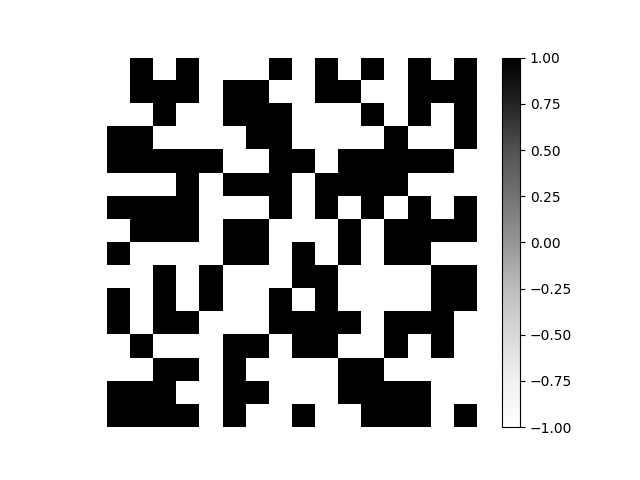}
  \caption{End state with no interaction}
\end{subfigure}
\caption{Sum of initial sentiment is $-4.031$ vs. sum of final sentiment after 1944 iterations is $106.869$. This indicates a dominantly positive (black) change in public sentiment.}
\label{fig-4}
\end{figure}

We implement a simple forward-Euler numerical scheme in {\em Python} to produce two numerical simulations that illustrate the change of individual initial sentiment after interacting with other members of a population in the presence of two attractive steady states.
These simulations appear in Figures \ref{fig-3} and \ref{fig-4}.
The simulation was run until the numerical solution approached an equilibrium.
In this case, each simulation ran until the maximum difference taken over each of the $16\times16$ squares between the last two iterations was less than $0.001$.
For comparison, we also provide an image of the steady-state as if no interaction (mixing) were present.
Polarization in public opinion is apparent in the development of pure white and pure black pools of color. 
Interestingly, we see the appearance of entire columns of black or white also develop.
This can be explained by the presence of sufficient influence in the interaction kernel.

\begin{figure}[htpb]
    \centering 
\begin{subfigure}{0.45\textwidth}
  \includegraphics[width=\linewidth]{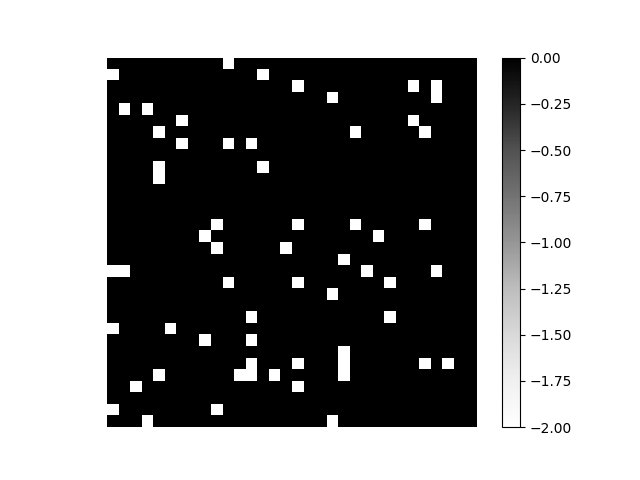}
  \caption{Emergent negative polarity: White pixels indicate black pixels changing to white.}
\end{subfigure}\hspace{1cm}
\begin{subfigure}{0.45\textwidth}
  \includegraphics[width=\linewidth]{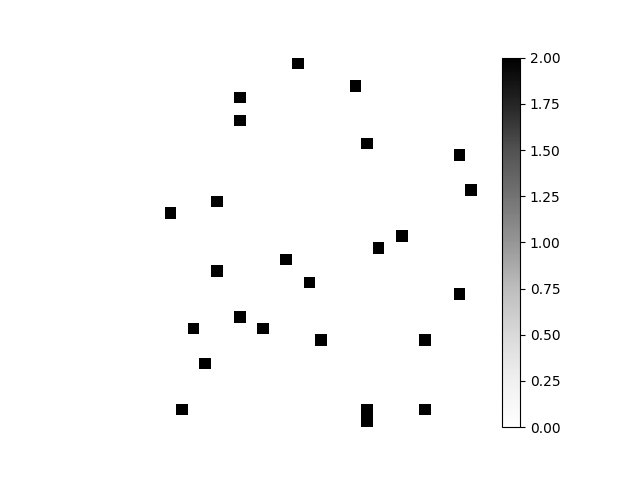}
  \caption{Emergent positive polarity: Black pixels indicate white pixels changing to black.}
\end{subfigure}\hspace{0cm}
\medskip
\begin{subfigure}{0.45\textwidth}
  \includegraphics[width=\linewidth]{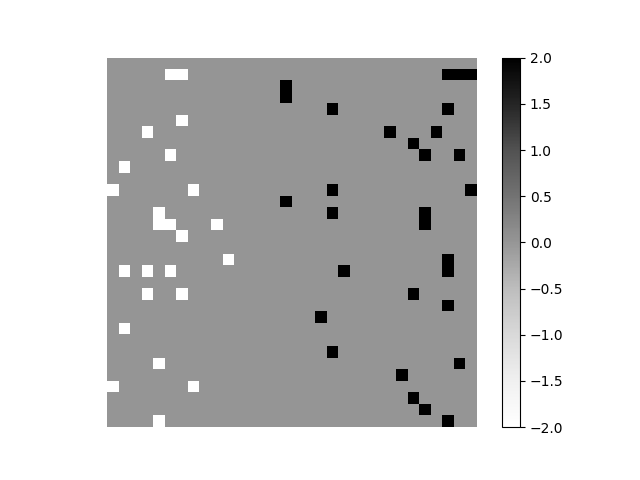}
  \caption{Emergent mixed polarity: Grouped/separated changes in polarity.}
\end{subfigure}\hspace{1cm}%
\begin{subfigure}{0.45\textwidth}
  \includegraphics[width=\linewidth]{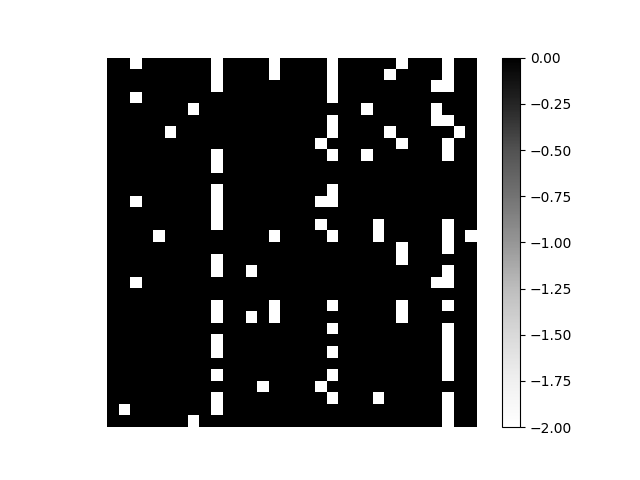}
  \caption{Emergent negative polarity: Strong changes with white bands in polarization.}
\end{subfigure}\hspace{0cm}
\medskip
\begin{subfigure}{0.45\textwidth}
  \includegraphics[width=\linewidth]{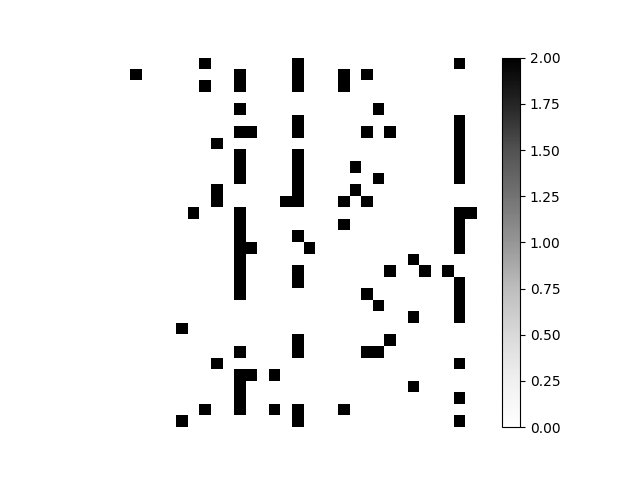}
  \caption{Emergent positive polarity: Strong changes with black bands in polarization.}
\end{subfigure}\hspace{1cm}
\begin{subfigure}{0.45\textwidth}
  \includegraphics[width=\linewidth]{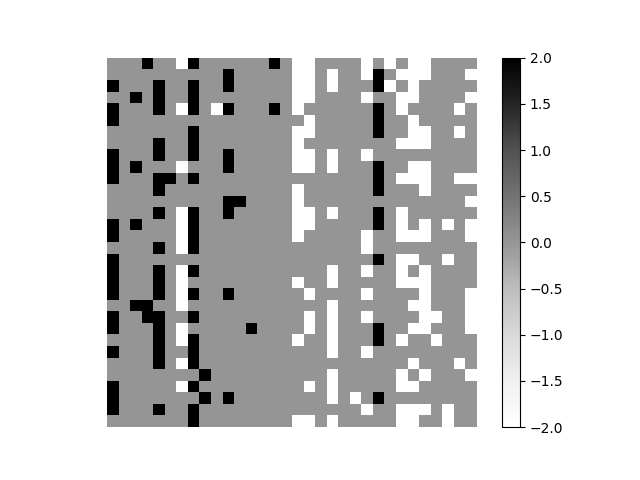}
  \caption{Emergent mixed polarity: Black and white bands in polarization.}
\end{subfigure}\hspace{0cm}%
\caption{Several $32\times32$ simulations produce these different difference maps.}
\label{fig-diff}
\end{figure}

Finally, we also report another method for detecting various degrees of sentiment change.
A so-called {\em difference map} records the difference between the initial condition and the equilibrium state. 
When a white pixel in an initial condition (which $=-1$ in our scheme) remains white in the equilibrium, the difference maps presents this pixel with the value $-1-(-1)=0$.
Similarly, when a black initial condition pixel remains black until the end, the difference map indicates $0$ for that pixel.
So, any pixel that does not change color will appear with a $0$ in the difference map.
Now we consider the case when change does occur.
This happens when a white pixel changes to black, or visa versa. 
The difference maps will show a white pixel that changed to black with the value of $1-(-1)=2$, and a black pixel that changes to white is valued with $-2$.
So far, subtle changes are difficult to detect with this method.
To help alleviate this, we define the difference map as follows,
\begin{equation}\label{eq:diff}
{\rm diff}(x) := {\rm sgn}\{\phi_N(x)\}-{\rm sgn}\{\phi_0(x)\}.
\end{equation}
Equation (\ref{eq:diff}) says the difference map value at pixel $x$ is the difference of the {\em sign} function applied to the equilibrium state (this gives values of only $-1$ or $+1$) and the sign function of the initial condition (again, only possible value here are $-1$ or $+1$).
The following difference maps illustrate a variety of effects produced by interactions modeled by nonlocal diffusion with random kernels.
These appear in Figure \ref{fig-diff}.
It should be noted that to better highlight the different emergent structures, we rely on a larger $32\times32$ system.

\section{Conclusions}

We have presented a nonlocal differential equation model to simulate the change and subsequent polarization in public sentiment.
The following observations and questions arise.
\begin{itemize}
\item {\em Dependency}: The columns of the sentiment data are implicitly assumed to represent independent questions. This should not always the case. The model could be fitted to express dependence between certain questions. This means we will also expect to see change along the (horizontal) rows during the sentiment evolution.
\item {\em Larger systems}: Thanks to modern computing architecture (GPUs) and modern programming libraries (Python's {\em NumPy}), the present algorithm is scalable to larger systems and even subject to more complicated interaction kernels.
\item {\em A sensitivity problem}: Which one sentiment change produces the greatest deviation from the current final sentiment? Solving this problem could help us understand and detect the propagation of misinformation. 
\item {\em Big data}: Beyond relying on people to fill out surveys, we can find other sources for our model's initial sentiment data. There is data collected by our smart phones, cellular servicer, computers, internet service provider, etc. 
\end{itemize}


\noindent{\bf Conflicts of Interest.} The author declares no conflict of interest. 

\noindent{\bf Funding.} This research did not receive any specific grant from funding agencies in the public, commercial, or not-for-profit sectors.

\bigskip

\bibliographystyle{plain}

\end{document}